\documentclass[a4paper,11pt]{amsart}
\usepackage{amsmath,amsthm,amsfonts,amssymb,amscd}
\usepackage{hyperref}
\usepackage{graphicx}

\theoremstyle{plain}
\newtheorem{Thm}{Theorem}[section]

\newtheorem{Lem}[Thm]{Lemma}
\newtheorem{Prop}[Thm]{Proposition}

\newtheorem{Rk}[Thm]{Remark}

\theoremstyle{definition}
\newtheorem{Defn}[Thm]{Definition}

\newtheorem{Ex}[Thm]{Example}

\newenvironment{pf}{ \begin{proof} }{ \end{proof} }

\DeclareMathAlphabet\EuScript{U}{eus}{m}{n}
\SetMathAlphabet\EuScript{bold}{U}{eus}{b}{n}

\DeclareFontFamily{U}{eus}{\skewchar\font'60}%
\DeclareFontShape{U}{eus}{m}{n}{<-6>eusm5<6-8>eusm7<8->eusm10}{}%
\DeclareFontShape{U}{eus}{b}{n}{<-6>eusb5<6-8>eusb7<8->eusb10}{}%

\newcommand{\Z}{\mathbb{Z}}

\newcommand{\R}{\mathbb{R}}
\newcommand{\C}{\mathbb{C}}

\DeclareMathOperator*{\SO}{\mathrm{SO}}

\DeclareMathOperator{\im}{im}
\DeclareMathOperator{\End}{\mathrm{End}}

\DeclareMathOperator{\aut}{\mathrm{Aut}}
\DeclareMathOperator{\Hom}{\mathrm{Hom}}

\DeclareMathOperator{\imag}{Im}
\DeclareMathOperator{\real}{Re}
\DeclareMathOperator{\torus}{T}
\DeclareMathOperator{\diff}{Diff}

\DeclareMathOperator{\spinc}{\mathrm{Spin}^c}

\newcommand{\ii}{\mathrm{i}}
\newcommand{\id}{\mathrm{id}}

\begin{document}

\title{Zero-sets of near-symplectic forms}
\author{Tim Perutz}
\address{DPMMS, Centre for Mathematical Sciences, University of Cambridge, Wilberforce Road, Cambridge CB3 0WB, U.K.}
\email{T.Perutz@dpmms.cam.ac.uk}
\date{Revised version: January 8 2007}
\begin{abstract} We give elementary proofs of two `folklore' assertions about near-symplectic forms on four-manifolds:
that any such form can be modified, by an evolutionary process taking place inside a finite set of balls, so as to have
any chosen positive number of zero-circles; and that, on a closed manifold, the number of zero-circles for which the
splitting of the normal bundle is trivial has the same parity as $1+b_1+b_2^+$.
\end{abstract}
\maketitle

\section{Statement of results}
In this paper we establish some properties of near-symplectic forms, as defined in \cite{ADK}:

\begin{Defn}
A two-form $\omega$ on an oriented four-manifold $X$ is called \emph{near-positive} if at each point $x\in X$,
either (i) $ (\omega\wedge \omega)(x) > 0$, or (ii) $\omega(x)=0$ and the intrinsic gradient
$(\nabla \omega )(x) \colon T_x X\to \Lambda^2 T_x^*X $ has rank 3 as a linear map.
It is \emph{near-symplectic} if, in addition, $d\omega=0$.
\end{Defn}

To clarify: at a point where $\omega(x)=0$, there is an intrinsic gradient
$(\nabla_v\omega)(x) \in \Lambda^2 T_x^*X$ in the direction $v$, for any $v\in T_x X$,
because $x$ is a zero of a smooth section of a vector bundle.
If $\omega\wedge \omega \geq 0$ in a neighbourhood of $x$ then
$(\nabla_v\omega)\wedge (\nabla_v\omega) \geq 0$ for all $v\in T_x X$.
The wedge-square quadratic form $\Lambda^2 T^*_x X\otimes \Lambda^2 T^*_x X \to \Lambda^4 T^*_x X \to \R^4 $, which is
defined up to a positive scalar, has signature $(3,3)$, so to say that $(\nabla \omega)(x)$ has rank 3 is to say that
its image is a maximal positive-definite subspace for the wedge-square form.

\begin{Lem}
The zero-set $Z_\omega
= \{ x \in X: \omega_x = 0 \}\subset X$ of a near-positive form $\omega \in
\Omega^2_X$ is a smooth $1$-dimensional submanifold.
\end{Lem}

\begin{pf}
Take $z\in Z_\omega$. Working over a small ball $B\owns z$, choose a three-plane subbundle $E\subset \Lambda^2 T^*B$
such that $E_z$ is complementary to $\im(\nabla\omega)(z)$.
Project $\omega|_B$ to a section $\overline{\omega}$ of $(\Lambda^2 T^*B)/E$.
Then $(\nabla \overline{\omega})(z)$ is surjective, so, shrinking $B$ if necessary, $\overline{\omega}$ vanishes transversely
along a 1-submanifold $\overline{Z}\subset B$. Moreover, $T_x \overline{Z}=\ker(\nabla\omega)(x)$
for $x\in \overline{Z}$, so $\omega$ is constant along $\overline{Z}$; hence $\overline{Z}=  Z_\omega\cap B$.
\end{pf}
A smooth path of two-forms $\{\omega_t\}_{t\in[0,1]}$ can be regarded as a section $\omega_\bullet$ of the pullback
of the vector bundle $\Lambda^2 T^*_x X$ to $[0,1] \times X$.
\begin{Defn}
Define a \emph{near-symplectic cobordism} on $X$ to be a
path $\{\omega_t\}_{t\in[0,1]}$ of closed two-forms such that for all $(x,t) \in X\times [0,1]$, either
(i) $ \omega_t(x)\wedge \omega_t(x)>0$, or (ii) $\omega_t(x)=0$, and
$(\nabla\omega_{\bullet})(t,x)$ has rank $3$. Here $(\nabla \omega_\bullet)(t,x)$ is the intrinsic gradient on
$[0,1]\times X$.
\end{Defn}
When $\omega_0$ and $\omega_1$ are near-symplectic, the zero-set
$Z= \bigcup_{t\in [0,1]}{ \{t\}\times Z_{\omega_t}}$ of the cobordism $\omega_\bullet$ is a smooth surface in
$[0,1] \times X$, realising a cobordism between $Z_{\omega_0}$ and $Z_{\omega_1}$.

In this article we prove the following:

\begin{Thm}\label{main}
Suppose that the connected, oriented four-manifold $X$ carries a near-symplectic form $\omega_0$ whose zero-set has $n$
components, where $n\geq 1$. Then
\begin{enumerate}
\item
there is a near-symplectic cobordism $\{ \omega_t\}_{t\in [0,1]}$, constant on the complement of a ball in $X$,
such that the zero-set of $\omega_1$ has $n+1$ components;
\item
if $n \geq 2$, there is a near-symplectic cobordism $\{ \omega_t\}_{t\in [0,1]}$, constant on the complement of a ball in
$X$, such that the zero-set of $\omega_1$ has $n-1$ components.
\end{enumerate}
In either case, the zero-surface $Z$ of the cobordism realises an
elementary cobordism between the zero-sets of $\omega_0$ and $\omega_1$:
the projection $Z \to [0,1]$ is a Morse function with just one critical point, of index 1.
\end{Thm}

The proof given here is elementary---it relies on stability arguments \emph{\`a la} Moser.
We shall show that given any pair of disjoint arcs embedded in the zero-set of $\omega_0$,
and a path in $X$ connecting interior points of these arcs, one can find a near-symplectic cobordism whose zero-surface
is the result of attaching a (2-dimensional) 1-handle to the pair of arcs so that the core of the handle is the chosen path.
Thus, given arcs on two distinct components of the zero-set, one can join them up to make a single component.
On the other hand, taking the arcs to to be on the same component,
the surgery increases by one the number of components.

A statement similar to Theorem \ref{main} is mentioned in Taubes' 1998 ICM lecture \cite{Ta1}, with
attribution to unpublished work of Karl Luttinger.
After I had written a draft of the present article, Taubes sent me a preprint, which has now appeared as \cite{Ta3},
in which he also shows that the number of zero-circles can be increased or decreased by a local modification.
His approach also uses 1-handle addition, but aside from that it is different.
It does not proceed through an evolution of the forms, and so does not include the cobordism statement,
but rather through a `mutation' in which one removes a ball whose intersection with the zero-set is a pair of arcs, then invokes Eliashberg's
classification of overtwisted contact structures to argue that it can be glued back so that the ends of the arcs
match up differently. Taubes also proves by means of a `birth' model that any near-symplectic form is
near-symplectically cobordant to one with one more zero-circle.

As a simpler variant, one can deduce the following weaker statement from a result of Calabi:
\begin{Prop}
Let $\omega$ be a near-symplectic form on an oriented four-manifold whose zero-set has $n\geq 0$ components.
Then there is a near-symplectic form whose zero-set has $n+2$ components and which coincides with $\omega$ outside an
embedded ball.
\end{Prop}
\begin{pf}
We consider near-symplectic forms $\omega_{f,g}$ on $\R^4=\R^3\times \R$ which arise from functions $f\in C^\infty(\R^3)$
and metrics $g$ on $\R^3$ which coincide with the Euclidean metric $g_0$ outside a compact set:
\[\omega_{f,g} =  d f\wedge dt  +  *_g df.  \]
The two-form $\omega_{f,g}$ is closed if and only if $df$ is $g$-coclosed. If $f$ is Morse then $\omega_{f,g}$ is near-positive,
and its zero-circles are in
bijection with the critical points. For example, when $f_0(x) = x_3$ we recover standard symplectic $\R^4$:
$\omega_{f_0, g_0} = \omega_{\R^4}$.

In \cite{Cal}, Calabi gives a recognition criterion for a closed one-form to be coclosed with respect to \emph{some} metric.
As an example, he shows that one can start with the function $f_0=x_3$ and modify it to a function $f_1$ which coincides with
$f_1$ outside a closed ball of chosen radius $r$, such that $df_1$ is $g_1$-coclosed for some $g_1$ standard outside the $r$-ball.

By Darboux's theorem, there exists $r>0$ such that the open ball $D^4(0;2r)$ of radius $2r$, with its symplectic form
$\omega_{f_0,g_0}$, embeds into our near-symplectic manifold. We can then plug in $\omega_{f_1,g_1}$ in place of
$\omega_{f_0,g_0}$, which has the effect of increasing by two the number of zero-circles.
\end{pf}

\begin{Ex}
When $X$  is closed, connected, and not negative-definite, there are near-symplectic forms with arbitrary positive
numbers of zero-circles. The existence of a near-symplectic form in this situation is well-known,
and can be proved either using Hodge theory
as in \cite{Leb,Ho2}---an argument which originated in the 1980s---or via handle decompositions and the flexibility of
overtwisted contact structures \cite{GK}.
\end{Ex}
\begin{Rk}
Let $\omega_0$, $\omega_1$ be two near-symplectic forms on a closed manifold $X$. Assuming either that $b_2^+(X)> 1$, or
that they represent the same cohomology class, these forms are near-symplectically cobordant (in $C^k$, $k<\infty$) by the
work of Honda \cite[Thm. 2.24]{Ho2}. When $\omega_0$
and $\omega_1$ are symplectic, the zero-set of a cobordism projects to an orientable surface in $X$ whose homology class is
presumably determined by $c_1(X,\omega_0)-c_1(X,\omega_1)$. Supposing this vanishes, one can ask whether there are
constraints on the topology of the surface, or on the function on it defined by the $t$-coordinate of the cobordism. For example,
when a three-manifold $M$ fibres over $S^1$, there is a standard construction of a symplectic form on $S^1 \times M$. Given two
possibly inequivalent fibrations on $M$, representing the same class in $H^1(M;\Z)$ and having the same
Euler class in $H^2(M;\Z)$, it is not clear whether the associated symplectic forms need be deformation-equivalent.
One could try first to use Morse theory on $M$ to construct a near-symplectic cobordism between them whose zero-set is a
union of product tori, and then to attempt to modify this cobordism so that its zero-set becomes empty.
\end{Rk}

I also learned the following topological statement from Taubes' lecture \cite{Ta1}; he attributes it to Gompf.
Though proofs are no doubt known to a few experts, no proof has appeared in print (and this one was worked out
independently).

A near-symplectic form $\omega$ determines a splitting of the normal bundle $N_{Z_{\omega}/X}$ into two sub-bundles,
on which the quadratic form $N_{Z_\omega/X}\to \underline{\R}$, $v\mapsto \langle \iota(z) \nabla_v \omega, v \rangle$
is respectively positive- or negative-definite; here $z$ is a non-vanishing tangent vector field on $Z_\omega$. One of these bundles
has rank 1. A component of $Z_\omega$ over which this line bundle is topologically trivial is called `even';
one for which it is non-trivial is called `odd'.

\begin{Thm}\label{Gompf}
The number of even circles in the zero-set of a near-symplectic form on a closed, connected four-manifold $X$ is
congruent modulo 2 to $1-b_1(X)+ b_2^+(X)$.
\end{Thm}
For example, there is an even number of them when $X$ admits an almost complex structure.

In the proof of Theorem \ref{main}, two chosen components are joined together
to form a single circle. The new circle has even parity if and only if the parities of the two old circles were different.

The zero-set $Z_\omega$ is always \emph{nullhomologous}: it is dual to the Euler class of $\Lambda^+_g$, which is
zero \cite{HH} (its components might not be nullhomologous, however).
Here is a brief indication of the proof of Theorem \ref{Gompf} in the
special case where the zero-circles are contained in disjoint balls $B_i$.
On the complement of the balls, the form distinguishes an almost complex structure $J$, up to homotopy.
This has a first Chern class, which extends uniquely to a class $c\in H^2(X;\Z)$.
The number $\frac{1}{4}(c^2-2e(X)-3\sigma(X))$, which reduces modulo 2 to $1-b_1+ b_2^+$, is the obstruction to
extending $J$ over the punctures. Each \emph{odd} circle makes an \emph{even} contribution to this obstruction, and vice versa.

The general case can be reduced to this one by carrying out surgery on all the circles so that they bound disjointly embedded discs.

\begin{Rk}
Given a near-symplectic form on $X$ with $m$ even and $n$ odd components, there is another with $m'$
even components and $n'$ odd components provided that $m-m'$ is even and $m'+n'>0$. To see this, first reduce the total number
of its zero-circles to $0$ or $1$. Pairs of even circles can then be introduced as in the proposition above. By Theorem
\ref{main}, a pair of even circles can be fused to give a single circle, which by Theorem \ref{Gompf} is necessarily odd.

The remaining question is to determine when the last remaining zero-circle (if it is odd) can be removed so as to obtain a
symplectic form.
\end{Rk}

\subsubsection*{Acknowledgements} Much of the work presented here appeared in the first chapter of my doctoral
thesis (Imperial College London, 2005), and I am indebted to my Ph.D. supervisor, Simon Donaldson, for useful discussions
about near-symplectic geometry. Thanks also to Carlos Simpson for sending me a copy of the unpublished manuscript \cite{LS};
it contains an example which we build on here.

\section{Local near-symplectic geometry}

\emph{A two-form is near-symplectic if and only if it is self-dual and harmonic for some conformal
structure and vanishes transversely as a section of the bundle $\Lambda^+$ of self-dual forms.} This follows from
the following lemma (see e.g. \cite{ADK}):

\begin{Lem}
(a) Let $X$ be a $4$-manifold with chosen conformal structure. Let $\omega\in \Omega^+_X$ be a transverse section of $\Lambda^+$.
Then $\omega$ is near-positive.

(b) Let $\omega$ be a near-positive form on an oriented $4$-manifold $X$. The following sets are in canonical bijection:
\begin{enumerate}
\item
conformal structures making $\omega$ self-dual;
\item
positive-definite three-plane sub-bundles $\Lambda^+\subset\Lambda^2_X $ having $\omega$ as section;
\item
almost complex structures $J$ on $X\setminus Z_\omega$ such that $g_J(\cdot,\cdot):=\omega(\cdot,J\cdot)$ defines a metric whose
conformal class extends to one on $X$.
\end{enumerate}

(c) These sets are non-empty, and contractible in the $C^\infty$ topologies.
\end{Lem}

The following self-dual form $\Theta$ on Euclidean $\R^4$ is the prototypical near-symplectic
form: let $(t,x_1,x_2,x_3)$ denote standard, positively-oriented coordinates.
Write
\begin{align}
&  \beta_1 = dt \wedge dx_1 + dx_2\wedge dx_3,  \notag  \\
&  \beta_2 = dt \wedge dx_2 + dx_3\wedge dx_1, \label{basic forms}\\
&  \beta_3 = dt \wedge dx_3 + dx_1\wedge dx_2. \notag
\end{align}
Define $\Theta  \in \Omega^2_{\R^4}$ by
\begin{equation}\label{Theta}
\Theta  = x_1\beta_1+ x_2 \beta_2 - 2 x_3 \beta_3,
\end{equation}
so that $\Theta^2 = 2(x_1^2+x_2^2+4 x_3^2)\, dt\wedge  dx_1\wedge  dx_2\wedge dx_3$. The form $\Theta$ is
near-symplectic with zero-set $Z_\Theta = \{x_1=x_2=x_3=0 \}$.

Two examples of near-symplectic forms with compact zero-set arise through symmetries of $\Theta$. These symmetries are
\[ \Theta(t, x)=\Theta(t-1,x) = \Theta(t, \sigma x ), \]
where $\sigma (x_1,x_2,x_3)=(x_1,-x_2,-x_3)$. By the first of these, $\Theta$ descends to a form
$\Theta_{\mathrm{ev}}$ on $S^1\times \R^3$. By the second, it descends to a form $\Theta_{\mathrm{odd}}$
on the mapping torus $\torus(\sigma) = (\R\times \R^3)/ (t-1,x)\sim (t,\sigma x) $.
Its zero-set is the zero-section of the mapping torus.

\subsection{Birth model (Luttinger and Simpson)} The following example is from \cite{LS} (apart from very slight variations
in the formulae): \emph{There is a path of closed forms $\{\omega_t\}_{t\in[-1/2,1/2]}$ on the unit ball in $\R^4$,
whose zero-set $Z_{\omega_t}$ is a circle for $t>0$, a point $\{p\}$ at $t=0$, and empty when $t<0$.
They are near-symplectic when $t\neq 0$. At $t=0$, $(\nabla\omega_0)(p)$ has rank $2$,
but $(\nabla^{\R^4\times \R}\omega_{\bullet})(p,0) $ has rank 3.
The path $\omega_\bullet$ is a near-symplectic cobordism.}

We write down such a path as follows. Using positively oriented coordinates $(x_1,x_2,x_3,x_4)$, and writing
$dx_{i j}$ for $dx_i\wedge dx_j$, define
\begin{equation} \label{zeta}
\zeta = x_2 (dx_{12}+dx_{34})  -  x_4 (dx_{14} + dx_{23}),
\end{equation}
\begin{equation}
\eta_t = \frac{1}{2}(x_1^2 + x_3^2 -t) (dx_{13} + dx_{42})  + x_4( x_1 dx_{12} +  x_3 dx_{32}).
\end{equation}

Fix any $\epsilon \leq 1/2$. Set $\omega_t= \zeta+ \epsilon\eta_t$. Then, on the open unit ball $D^4 \subset \R^4$, the
path $\{\omega_t\}_{t\in[-1/2,1/2]}$ has the properties just described.
To see this, observe that $\omega_t$ is a $C^\infty(\R^4)$-linear combination of the three forms
\[  dx_{12}+dx_{34}, \quad dx_{14}+dx_{23}-\epsilon (x_1 dx_{12}+x_3 dx_{32}) ,\quad  dx_{13}+dx_{42}.\]
A simple determinant computation shows that these span a maximal positive-definite subbundle of the two-forms over $D^4$
when $\epsilon\leq 1/2$. Hence $\omega_t$ has positive square except at those points where the three coefficients
(namely $x_2$, $x_4$ and $(x_1^2+x_3^2-t)$) are zero.
At those points, $\partial_2{\omega_t}$ and $\partial_4 {\omega_t}$ are linearly independent,
so $\nabla \omega_t$ has rank $\geq 2$. By including also $(x_1\partial_1 + x_3\partial_3)\omega_t$
one finds that $\nabla \omega_t$ has rank 3 along the zero-set, except when $t=0$ and $x=0$.
In that case, the $t$-derivative is non-zero.

\subsection{Surgery model} We vary the birth model by replacing the family of forms $\eta_t$ by a new family $\eta_t'$:
\begin{equation} \label{surgery model}
\eta_t' = \frac{1}{2}(x_1^2 - x_3^2 -t) (dx_{13}+ dx_{42}) +  x_4( x_1 dx_{12}  - x_3dx_{32}).
\end{equation}
Taking $ t\in[-1/2,1/2]$ and $\epsilon\leq 1/2$, the forms $\zeta+\epsilon \eta_t'$ are near-symplectic on $D^4$
(except when $t=0$, where the rank condition fails at the origin).
However, the zero-set is now $\{ (x_1,0,x_3,0): x_1^2 - x_3^2 -t =0\}$. Thus as $t$ passes through zero, the zero-set undergoes a
surgery (it changes by addition of a one-handle). When $t=0$, $\nabla (\zeta+\epsilon \eta_0')(p)$ has rank $2$ at the
double point $p$ of the zero-set. This model is again a near-symplectic cobordism.

\subsection{Geometry near the zero-circles}\label{near circles}

Taubes \cite{Ta2} has established the following picture of the differential geometry of a near-symplectic form $\omega$ near its zero-set $Z=Z_\omega$. Fix a metric $g$ with $*_g\omega=\omega$, and identify the intrinsic normal bundle $N_{Z/X}$ with the metric complement $(TZ)^\perp$.

{\bf (a)} $\nabla\omega$ defines a vector bundle isomorphism $N_{Z/X} \to\Lambda^+|_Z$.
In particular, there is a universal procedure for orienting $Z$ (as for the zero-manifold of any transverse section of an
oriented vector bundle).

{\bf (b)} Let $z$ (resp. $\zeta$) denote the unit-length oriented vector field (resp. one-form) on $Z$. The map
\[\Lambda^+|_Z \to (N_{Z/X})^*,\quad \eta\mapsto \iota(z)\eta\]
is a bundle isomorphism; its inverse is $\alpha\mapsto (\zeta \wedge\alpha)+*(\zeta\wedge\alpha)$.

{\bf (c)} We have a sequence
\[\begin{CD}
N_{Z/X}         @>{\nabla\omega}>> \Lambda^+|_Z
                @>{\iota(z)}>> (N_{Z/X})^*
                @>{\text{metric}}>> N_{Z/X}
\end{CD}, \]
whose composite we denote by $S_{\omega,g}\in \End(N_{Z/X})$. The endomorphism $S_{\omega,g}$ is a
\emph{self-adjoint, trace-free automorphism}.

{\bf (d)} It follows that, at each point of $Z$, $S_{\omega,g}$ has a basis of eigenvectors, and that two of its
eigenvalues are positive and one negative (this statement pins down our convention for the universal orientation rule).
The positive and negative eigenspaces trace out eigenbundles $L^+_g$ and $L^-_g$
(here $L^-_g$ has rank 1 and is obviously locally trivial; $L^+_g$ is its orthogonal complement).
Thus $N_{Z/X} = L^+_g\oplus L^-_g$.

The auxiliary choice of metric is actually irrelevant. Whilst there are many conformal classes $[g]$ which
make $\omega$ self-dual, they are all the same along $Z$, because for $z\in Z$, they
satisfy $\im(\nabla \omega)(z)=\Lambda^+_{[g],z}$.
The sub-bundle $TZ^\perp\subset TX|_Z$ is therefore intrinsic to $\omega$.

We condense these points into a metric-free statement:

\begin{Prop}
A near-symplectic form $\omega$ determines
\begin{enumerate}
\item
a canonical orientation for the 1-manifold $Z_\omega$;
\item
a canonical embedding of the intrinsic normal bundle $N_{Z_\omega/X}$ as a sub-bundle of $TX|_Z$ complementary to $TZ$;
\item
a conformal class of quadratic forms $S_\omega$ on $N_{Z_\omega/X}$, of signature $(2,1)$;
\item
a vector bundle splitting $N_{Z_\omega/X} = L^+\oplus L^- $, such that $S_\omega$ is positive-definite on $L^+$
and negative-definite on $L^-$.
\end{enumerate}
\end{Prop}
A trivial three-plane bundle over $S^1$ admits precisely two topological splittings, distinguished by the orientability
of the summands. Thus a basic invariant of a near-symplectic form $\omega$ is its `parity function'
$  \varepsilon_\omega: \pi_0(Z_\omega) \to \Z/2$, which sends a component-circle $\Gamma$ of $Z_\omega$
to $\langle w_1(L^-), [\Gamma] \rangle.$

This function is equivariant under oriented diffeomorphisms, and invariant under isotopies
of near-symplectic forms $\omega_t$ (note that the topology of the zero-set cannot change if \emph{all}
the $\omega_t$ are near-symplectic).

Each of the model forms $\Theta_{\mathrm{ev}}$ and $\Theta_{\mathrm{odd}}$ defined afer (\ref{Theta}) vanishes
along a circle, and their parities are those suggested by the notation.

\section{Moser-Honda deformation arguments}

The extent to which the Moser deformation argument from symplectic geometry applies in the near-symplectic setting was worked out
by Ko Honda \cite{Hon}. He showed that any near-symplectic form can be deformed, through near-symplectic forms with fixed zero-set,
so that near a chosen zero-circle it embeds into $\Theta_{\mathrm{ev}}$ or $\Theta_{\mathrm{odd}}$.\footnote{One can do
without the preliminary deformation at the price of a loss of smoothness. Given two near-symplectic forms $\omega_0 $, $\omega_1$,
with zero-circles $Z_0$, $Z_1$ of the same parity, one has $\psi^*\omega_1=\omega_0$ near $Z_0$ for some map
$\psi\colon \mathrm{nd}(Z_0)\to \mathrm{nd}(Z_1)$; however, $\psi$ is in general only Lipschitz continuous along $Z_0$ (it is $C^\infty$ elsewhere).}
We use Honda's method---a careful application of Moser's argument on the complement for the zero-set---to prove a slight
variant of this result:

\begin{Lem}\label{Honda}
Let $\omega\in \Omega^2_X$ be near-symplectic, and let $A$ be a closed arc embedded in its zero-set $Z_\omega$.
Let $U$ be a neighbourhood of $A$. There exists a smooth family $\{\omega_s\}_{s\in[0,1]}$ of near-symplectic forms,
with constant zero-set, such that $\omega_s (p) = \omega (p)$ if $s=0$ or $p \in X\setminus U$, and where $\omega_1$ is standard near
$A$ in the following sense: there is an orientation-preserving open embedding $i\colon [-1,1]\times D^3(0;r)\to X $,
for some $r>0$, mapping $[-1,1]\times\{0 \}$ to $A$, such that $i^* \omega_1 = \Theta$.
\end{Lem}
Here $D^3(0;r)$ denotes the open ball of radius $r$ centred at $0\in \R^3$.
\begin{pf}

{\bf 1.} There is no loss of generality in assuming that $X= (-2,2)\times D^3(0;r)$;
that the zero-set of $\omega$ is $Z= (-2,2)\times\{0\}$; that $A = [-1,1]\times\{0\}$;
and that $U = (-1-\epsilon,1+\epsilon)\times D^3(0;r')$, for some $r' \leq r$. Writing $t$
for the coordinate on $(-2,2)$, we may suppose that $\partial_t$ is a positively oriented
vector field on $Z$ for the orientation determined by $\omega$.

Take a metric $g$ for which $\omega$ is self-dual. We can then find an orthonormal frame $(e_1,e_2,e_3)$
for $N_{Z/X}$ such that $L^+_g=\mathrm{span}(e_1,e_2)$ and $L^-_g = \mathrm{span}(e_3)$.

The metric and the choice of $e_i$ give rise to normal coordinates $(t,x_1,x_2,x_3)$ near $Z$,
identifying a small neighbourhood with $ (-\epsilon,1+\epsilon) \times D^3(0;r'')$. We then have
three basic two-forms $\beta_i$, as in (\ref{basic forms}). In the basis $(e_1,e_2,e_3)$,
the linear map $S_{\omega,g}(t)$ (see Section \ref{near circles}) is represented by a trace-free symmetric matrix
$S(t)=S^+(t)\oplus S^-(t)$, where $S^+(t)$ is $2\times 2$ and positive-definite, and $S^-(t)<0$.
Regarding $x$ and $\beta$ as column vectors, we have an expansion near $Z$:
\begin{align*}
\omega(t,x)     &= x\cdot S(t)   \beta +O(|x|^2)  \\
                        &= (x_1 \: x_2) \, S^+(t) \left(\begin{array}{c}
                                                                        \beta_1\\
                                                                        \beta_2
                                                                        \end{array}\right)
                                                        + x_3 S^-(t) \beta_3 + O(|x|^2).
\end{align*}
The advantage of paying attention to the eigenbundle splitting is that one can then apply convexity arguments.

Define $\theta_s  = (1-s)\omega+s(x_1\beta_1 + x_2\beta_2 - 2x_3\beta_3)$.
On a small neighbourhood of $Z$, this is a family of near-symplectic forms, with common zero-set $Z$.

{\bf 2.} On the total space of a vector bundle with a chosen metric one has the following explicit Poincar\'e
lemma. Let $p \colon V \to M$ be the bundle projection, $i: M\to V$ the zero-section.
Let $h_t: V\to V$ be the fibrewise dilation $x \mapsto tx$, and $R$ the Euler vector field, that is, the vertical
vector field defined on any fibre $V_x$ by $R(y)=\sum{y_j \frac{\partial}{\partial y_j}}$, where the $y_j$ are orthonormal Euclidean
coordinates on $V_x$. The map
\[ H\colon  \Omega^*_V\to \Omega^{*-1}_V;\quad \alpha\mapsto \int_0^1{h_t^*(\iota(R)\alpha) t^{-1} dt}  \]
satisfies $\alpha-p^*i^*\alpha=d (H\alpha)+H d\alpha$.
Notice that if $\alpha$ vanishes along $i(M)$ then $H(\alpha)$ vanishes to second order along $i(M)$.

Applying this to a neighbourhood of the zero-section in $N_{Z/X}$, we find that $\eta_s := H (\theta_s )$ vanishes to second order
along $Z$, and satisfies $d\eta_s =\theta_s$.

{\bf 3.} Let $U'\subset U$ be a neighbourhood of $A$ on which $\theta_s$ is defined and near-symplectic.
Take a non-negative cutoff function $\chi$, supported in $U'$ and identically 1 near $A$.
We choose $\chi$ to have the shape $\chi(t,x) = \chi_1(t)\chi_2(|x|)$ for functions $\chi_1, \chi_2\in C^\infty(\R)$.
Introduce the vector fields $v_s$ on $X\setminus Z$ defined by $\iota(v_s)\theta_s + \chi \eta_s =0 $.

The family $\{v_s\}_{s\in [0,1]}$ generates a flow $\{\phi_s\}_{s\in [0,1]}$ on $X\setminus Z$, supported inside $U$.
The reason is that, along $Z$, $\nabla \dot \theta_s (u)$ is non-zero for all $0\neq u\in N_{Z/X}$, while $\nabla (\chi\eta_s) (u)=0$, so that $|v_s(x)  | \leq k |x|$ near $Z$, for some constant $k$. A trajectory $x_s$, defined on some
interval $[0,s']$, satisfies $\frac{d}{ds}(\log |x_s|)\geq -k$; integrating over $[0,s']$, we obtain $|x_{s'}|\geq e^{-ks'} |x_0|$.
This shows that the trajectory stays inside $X\setminus Z$, and therefore that $\phi_s$ is defined.

Now define $\omega_s'$ by setting $\omega_s' = \phi_s^* \omega$ on $X\setminus Z$ and $\omega_s'(z)=0$ for $z\in Z$.
Outside $U$ we have $\omega_s'=\omega$. Moser's argument shows that $\omega_s'=\theta_s$ in some neighbourhood of $A$.

The path $\omega_s'$ gives an isotopy with all the required properties, except one: $\omega_s'$
fails to be smooth at certain points of $Z$---the points where $\chi_1(t)$ is non-constant.
At those points, its coefficients might merely be Lipschitz.

{\bf 4.} The `problem-points' of $\omega_s'$ lie strictly inside a pair of regions
$R_- = I_-\times \bar{D}^3(0;r_1)\subset U \setminus A$ and $R_+=I_+\times \bar{D}^3(0;r_2)\subset U \setminus A$,
for closed intervals $I_-\subset (-1-\epsilon,-1)$ and $I_+ \subset (1,1+\epsilon)$.
The map $\phi_s$ embeds $R_-\cup R_+$ into $X$, and acts as the identity along $Z$.
We can extend these embeddings $\phi_s$ to a path of homeomorphisms $s\mapsto \psi_s\in \mathrm{Homeo}(X)$,
such that (i) $\psi_0 = \id$; (ii) $\mathrm{supp}(\psi_s)\subset U$; (iii) $\psi_s|_{(R_-\cup R_+)} =\phi_s|_{(R_-\cup R_+)}$;
(iv) the only points where smoothness fails (smoothness of $\psi_s$, or of the path) are inside $R_+\cup R_- $;
and (v) $\phi_s$ is the identity in a neighbourhood of $A$.

Finally, setting $\omega_s = \psi_s^{-1*}\omega_s'$, we obtain a smooth path with the required properties.
\end{pf}

\section{Proof of Theorem \ref{main}}

Let $\theta$ be near-symplectic, and $Z_0$, $Z_1$ two components of $Z=Z_\theta$, not necessarily distinct.
Take distinct points $p_0\in Z_0 $ and $p_1\in Z_1$, and a smooth
path $\gamma\colon [0,1]\to X$ from $p_0$ to $p_1$, transverse to $Z$, with $\gamma^{-1}(Z)=\{0,1\}$.
The plan is this: first, we show how to modify $\theta$ so that it is standard in some neighbourhood $N$ of $\im(\gamma)$.
We then exhibit a path of forms on $N$,
constant near $\partial N$, which realises a surgery on zero-sets; these can be plugged into $\theta$.

The delicacy lies in finding a variant of our earlier `surgery model' which is a constant path of forms outside a
suitable neighbourhood of the origin. This is awkward partly because we cannot control the shape of the neighbourhood $N$:
in local coordinates we can take it to be an ellipsoid, but it might be very thin.

\begin{Defn}
A \emph{letter H} in a near-symplectic manifold $(X,\theta)$ consists of three smooth embeddings
$v_1,v_2,h\colon [0,1]\to X$, such that $v_1(1/2)= h(0)$ and $v_2(1/2)=h(1)$ (transverse intersections),
with no other intersection points. We require that $\im(v_1)$ and $\im(v_2)$ lie in the zero-set $Z_\theta$.
The interior of the arc $\im(h)$ joining them then necessarily consists of symplectic points.
\end{Defn}
\begin{Lem}
Let $H_i$ be a letter H in $(X_i,\theta_i)$, for $i=0$, $1$. Then $\theta_0$ is isotopic, through near-symplectic forms with
fixed zero-set, to a form $\theta_0'$, such that there are regular neighbourhoods $\mathrm{nd}(H_i) \supset H_i$ and a
diffeomorphism $\phi\colon \mathrm{nd}(H_0) \to \mathrm{nd}(H_1) $ with $\phi^* \theta_1=\theta_0'$.
\end{Lem}
\begin{pf}
We may assume that we are in $\R^4$, and that $H_1$ coincides with $H_2$.
By Lemma \ref{Honda}, we may also assume that $\theta_1$ and $\theta_2$ coincide in a neighbourhood of the two vertical arcs.

Let $P$ be the `horizontal stroke' $\im(h)$ in the letter H. For definiteness, say $P =  [-1,1]\times \{ (0,0,0)\} \subset \R^4$.
We claim that there exist linearly independent vector fields $(v_1=\partial_{x_1},v_2,v_3,v_4)$
along $P$ such that the matrix $(\theta_0(v_i,v_j))_{ij}$ takes the form
\[ \lambda(t)\left(
\begin{array}{cc}
  J    &  0  \\
  0   &   J   \\
\end{array}
\right),
\quad J=\left( \begin{array}{cc} 0 & -1\\ 1& 0\end{array}\right), \quad
\lambda(t)\geq 0; \]
and also linearly independent vector fields $(v_1'=\partial_{x_1},v_2',v_3',v_4')$ along $P$ such that
the matrix $(\theta_1(v_i',v_j'))$ takes the same form (for the same function $\lambda$).
Moreover, we may assume that, near $\partial P$, where $\theta_0$ and $\theta_1$ coincide,
we have $(v_1',v_2',v_3',v_4')= (v_1,v_2,v_3,v_4)$.

To see this, decompose $P$ into three subintervals:
$P_1=[0,\delta ] \times \{(0,0,0)\}$, $P_2 = ([\delta,1-\delta])\times \{(0,0,0)\}$, $P_3 = [1-\delta,1]\times \{(0,0,0)\}$.
If $\delta$ is small, we can find such frames over $P_1\cup P_3$, and we may assume, rescaling as necessary,
that the function $\lambda$ is identically equal to $1$ near the interior endpoints of $P_1$ and $P_3$. Over $P_2$,
we can find symplectic trivialisations of the tangent bundle, for either form, and since the symplectic group is path-connected,
we can make these trivialisations match those on $P_1$ and $P_3$.

There is now an open neighbourhood $N_1 \supset P$ and a diffeomorphism $\Phi\colon N_1\to N_1$, such that
$\Phi (x)=x$ if $x\in P$ or $x$ lies close enough to $\partial P$, and such that, along $P$, $D \Phi (v_i)=v_i'$. Define
\[ \theta_s = (1-s)\theta_0 + s \Phi^* \theta_1 \in \Omega^2_{N_1},\quad s\in [0,1]. \]
There is a smaller open neighbourhood $N_2\supset P$ on which these forms are near-symplectic.

The function $s\mapsto \theta_s\in \Omega^2_{N'}$ has constant derivative $d\theta_s /ds = \Phi^*\omega_0 - \theta_0$.
We may write $d\theta_s / ds = d\alpha$, where $\alpha\in \Omega^1_{N_2}$ vanishes near $\partial P$.
Since $\theta_s$ is a constant path along $P$, we may assume further that $\alpha$ vanishes along $P$.

Take an open set $N_3\supset P$ with $\overline{N_3} \subset N_2$, and a cutoff function $\R^4\to [0,1]$
such that $\chi(N_3)=\{1\}$ and $\chi(\R^4 \setminus  N_2)=\{0\}$. Define a family of vector fields $v_s$ on $\R^4$,
supported on $N_3$, by $\iota (v_s) \theta_s +\chi \alpha =0$. These fields are zero near the common zero-set of the $\theta_s$.
The family $\{v_s\}$ generates a flow $s\mapsto F_s$, which acts as the identity along $P$.
The near-symplectic form $F_1^*\theta_0$ coincides with $\theta_0$ outside $N_2$
and with $\theta_1= \Phi^*\theta_0$ in some still smaller neighbourhood of $P$, say $N_4 \subset N_3$.
The result follows.

\end{pf}

\subsection{Sharpening the surgery model}

The forms appearing in our surgery model fit into a larger family, which contains different examples more
amenable to embedding in other manifolds.

For any smooth function $f\colon \R^2\to \R$, set
\[ \eta^f = f(x_1,x_3)  (dx_{13}+dx_{42})+ x_4 (\frac{\partial f}{\partial x_1}dx_{12}+ \frac{\partial f}{\partial x_3}dx_{32});\]
\[ \omega^f = \zeta+ \frac{1}{2}\eta^f. \]
Here $\zeta$ is as in (\ref{zeta}). Note that $\eta^f$ and $\zeta$ are closed. If the condition
\begin{equation}\label{positivity} \left | x_4^2 \frac{\partial f}{\partial x_3} + x_2 x_4 \frac{\partial f}{\partial x_1}\right|\leq x_2^2 +x_4^2
 \end{equation}
holds then $\omega^f\wedge \omega^f  \geq (x_2^2+x_4^2)dx_{1234}$. Moreover,
$(\omega^f\wedge \omega^f)(x_1,0,x_3,0) = 2 f(x_1,x_3)^2 dx_{1234}$.
Hence $\omega^f$ satisfies the dichotomy that at any point it either has positive square or is zero.

We now work in $ U_\delta = \{ |x_1|<  2, |x_3|<2, x_2^2+x_4^2 < \delta \}$, where $\delta\leq 1$.
Put
\begin{equation}\label{function f}
f(x_1,x_3) =\frac{1}{4}(x_3^2-1).
\end{equation}
For this $f$, condition (\ref{positivity}) holds, and $\omega^f$ is near-symplectic in $U_\delta$.
Its zero-set consists of the two straight line segments, $\{ x_3=\pm 1, x_2=x_4=0\}$.

\begin{Lem}
There is a smooth family of functions $f_t^\delta\colon [-2,2]\times[-2,2]\to \R$, parametrised by $\delta\in
(0,\frac{1}{10}]$ and $t\in[-1,1]$, with the following properties:
\begin{itemize}
\item
$f_{-1}^\delta = f$, the function appearing in (\ref{function f});
\item
$f_t^\delta$ is independent of $t$ outside the square $[-1-2\delta, 1+2\delta]\times [-1,1]$;
\item
when $t< 0$, the zero-set of $f^\delta_t$ consists of two arcs, homotopic to those of $f$; when $t=0$, the zero-set is a letter X; and when $t>0$ it again consists of two arcs, connecting the endpoints the other way. Moreover, the vanishing is of order exactly 1, except at the singular point when $t=0$ where it is quadratic; and
\item
$f^\delta_t$ is bounded in $C^1$, uniformly in $t$ and $\delta$.
\end{itemize}
\end{Lem}
\begin{pf}
See Figure 1 for contour plots showing what such a family might look like. They show five different $t$-values for a fixed
$\delta$. Each box represent the square $[-1-2\delta,1+2\delta]\times [-1-2\delta,1+2\delta]$.
The functions are $t$-independent outside this range, and $f_{-1}^\delta=f$, so the values of the functions
on the contours are implied by the points at which the contours hit the boundary of the square (all contours do hit the
boundary). The bold contours are the zero-sets.

Begin by constructing smooth functions $f_1^\delta\colon [-2,2] \to \R$ as in the final picture, taking the zero-sets
to be $\delta$-independent and the functions to be even in both variables.
Then $\sup | f_1^\delta | = 3/4$ and the derivative of $f_1^\delta$ is uniformly bounded
(note that it is at most $O(\delta)$ on the region where $f_1^\delta\leq 0$). Now define
$f_t^\delta= \frac{1-t}{2} f + \frac{1+t}{2}f_1^\delta$. It is not immediately clear that the contours of $f^\delta_t$ resemble those in the
pictures. However, if we choose $f_1^\delta$ so that $\partial f_1^\delta/\partial x_1$ vanishes only where $x_1=0$,
we find that for $t>-1$, the only zero of $f_t^\delta$ in $(-1-2\delta,1+2\delta)\times (-1-2\delta,1+2\delta) $
which is a critical point of $f$ is the origin. Moreover, the origin is critical for
just one $t$-value (shown as $t=0.5$ in Figure 1), and this is a non-degenerate critical point. For all other
$t$-values, the zero-sets are smooth, and hence evolve by smooth homotopy. It follows that the family $\{f_t^\delta\}$
satisfies the requirements of the lemma.
\end{pf}

\begin{centering}
\begin{figure}[t!]
\includegraphics[width=13cm]{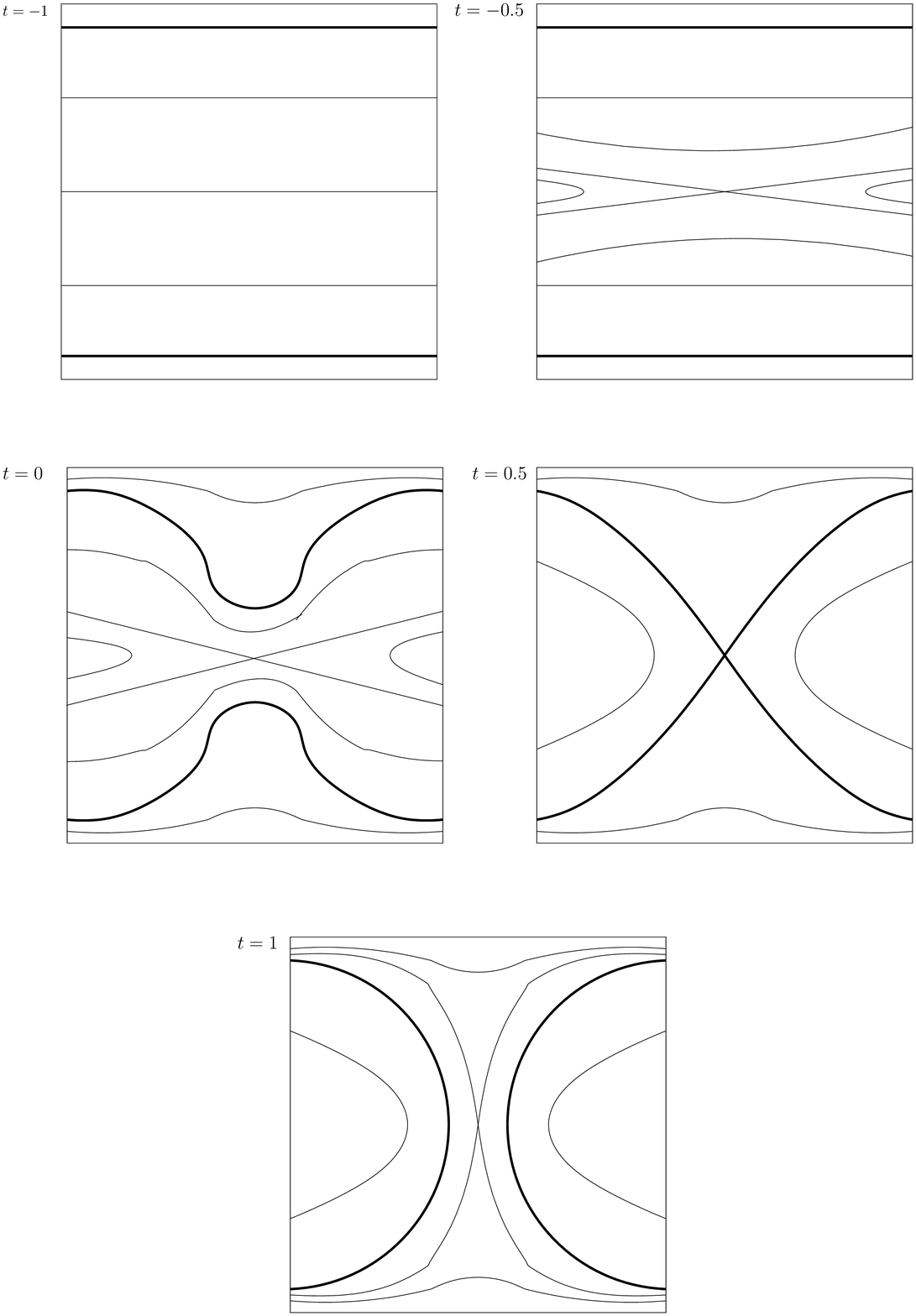}
\caption{\emph{Contour plots for the functions $f^\delta_t$, as $t$ varies from $-1$ to $1$.}}
\end{figure}
\end{centering}
Returning to the construction of our refined surgery model, we fix $\{f_t^\delta\}$, and determine an $\epsilon\in (0,1]$
such that the functions $\epsilon f^t_\delta$ all satisfy condition (\ref{positivity}) in $U_\delta$.
This is possible because of the uniform $C^1$ bound.

Given a letter H in a near-symplectic manifold $(X,\theta)$, we can perturb $\theta$ with fixed zero-set to a new form
$\theta'$, and embed the form $\omega^{\epsilon f}$ in the region
$ \{ |x_1|\leq \delta, |x_3|\leq 1+\delta,  x_2^2+x_4^2 < \delta \}$ into $(X,\theta')$ along the H.

We choose to reparametrise $\R^4$ by the dilation $D_\delta\colon (x_1,x_2,x_3,x_4)\mapsto (\delta x_1,x_2,x_3,x_4)$,
which leaves $f$ invariant;
so we have an embedding of $D_\delta^* \omega^{\epsilon f} $ on the region $U_\delta$ into $(X,\theta)$. We have
\[ D_\delta^* \zeta =  x_2(\delta dx_{12}+dx_{34}) -x_4(\delta dx_{14}+dx_{23}),\]
\[  D_\delta^* \eta^f = \frac{1}{4}(x_3^2-1) (\delta dx_{13}+dx_{42}) + \frac{1}{2}x_3x_4 dx_{32}.  \]
Consider the family of forms
\begin{equation}\label{final model}
 \Omega^\delta_t: = D_\delta^* \zeta +  \frac{\delta \epsilon}{2} \left(
	f_t^\delta(x_1,x_3)  (\delta dx_{13}+dx_{42})+ x_4 (\frac{\partial f_t^\delta}{\partial x_1}dx_{12}
	+ \frac{\partial f_t^\delta}{\partial x_3}dx_{32})\right).
\end{equation}
We observe, for any fixed $\delta$, the path $t \mapsto \Omega^\delta_t$ is a near-symplectic cobordism.
Indeed, the estimate (\ref{positivity}) for  $\epsilon f^t_\delta$ implies that
$( \Omega^\delta_t)^2 \geq \mathrm \delta (x_2^2+x_4^2)dx_{1234}$, while on the plane where $x_2=x_4=0$ we have
$( \Omega^\delta_t)^2 = \frac{\delta^3 \epsilon ^2}{2} (f^\delta_t)^2$.

\subsection{Stability arguments}
The formula (\ref{final model}) defines a path of forms $\Omega_t^\delta$, independent of $t$ when $(x_1, x_3)$ lies outside a certain
neighbourhood of $0\in\R^2$. Moreover, $\Omega_{-1}^\delta$ embeds in our near-symplectic manifold, for some $\delta$.
We would now like to refine the path so as to make it constant outside a suitable neighbourhood of the origin in $\R^4$.
We will achieve this using a stability principle from contact geometry---effectively, a refined Moser argument:

\begin{Prop}
Let $\{ \omega_t \}_{t\in [0,1]}$ be a smooth path of symplectic forms on a manifold $M$, and let $H\subset M$ be a compact, orientable hypersurface. Suppose that there exists a smooth path of one-forms $\alpha_t$ on $M$, with $d\alpha_t =\omega_t$, all of which are contact forms on $H$. Then there exists a neighbourhood $N$ of $H$, and smoothly varying embeddings $\phi_t\colon N\to M$, such that (i) $\phi_0$ is the inclusion; (ii) $\phi_t(H)=H$; and (iii) $\phi^*\omega_t=\omega_0$.
\end{Prop}
\begin{Rk}\label{rephrasing}
By carrying out a topologically trivial surgery on a tubular neighbourhood of $H$, we can rephrase the conclusion as follows: there is a new path $\omega'_t$, with $\omega_0'=\omega_0$, which coincides with the old one outside a small neighbourhood of $H$, and which is constant on $H$.
\end{Rk}
\begin{pf}
Let $\lambda_t$ be the Liouville vector field, characterised by $\iota(\lambda_t)\omega_t=\alpha_t$. It defines an embedding
\[ e_t\colon  [-\epsilon, \epsilon] \times H \to M,\quad (s, x)\mapsto \mathrm{flow}_{s \lambda_t} (x),  \]
where $\epsilon$ is chosen small enough that the flow of $\lambda_t$ exists on $[-\epsilon,\epsilon]$, for all $t$, and we have
$e_t ^* d( \mathrm{e}^s \alpha_t) = \omega_t$. By Gray's stability theorem, there is a path $t\mapsto f_t \in \diff(H)$,
with $f_0=\id$, satisfying $ f_t^*\alpha_t = \mathrm{e}^{g_t} \alpha_0 $ for some functions $g_t$. Put
\[ \phi_t =  e_t \circ (\mathrm{e}^{-g_t} \times f_t) \circ e_0^{-1}  \]
on a neighbourhood $N$ of $H$ where $e_0^{-1}$ is defined. Then $\phi_t$ restricts to $H$ as $f_t$, and $\phi_t^*\omega_t =\omega_0 $.
\end{pf}
We note that if the path $\alpha_t$ is constant in a (Liouville) tubular neighbourhood $N(H')$ of a region $H'\subset H $, then the
resulting embeddings $\phi_t$ are inclusions in $N(H')$.

\begin{centering}
\begin{figure}
\includegraphics[width=13cm]{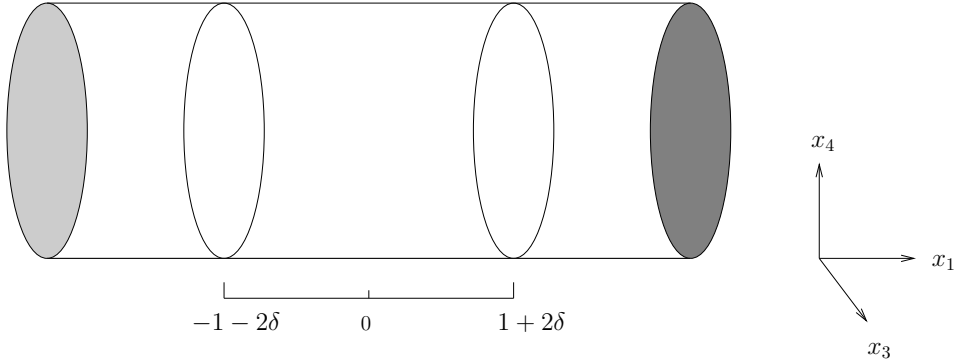}
\caption{\emph{Diagram of an $(x_1,x_2,x_4)$-section. The stability argument is to be applied to the
cylindrical hypersurface $C_\delta$ depicted here. The path of forms is constant except in the portion $C'_\delta$ between
the inner ellipses.}}
\end{figure}\end{centering}

We shall apply this to (a part of) the boundary of the region $U_\delta$.
The portion of $\partial U_\delta$ of interest is the subset
\[ C_\delta = \{|x_1| \leq 2, \; |x_3|\leq 2, \; x_2^2 +x_4^2 =\delta^2 \}. \]
We consider the path of forms $\Omega^\delta_t$ near $\partial U_\delta$.
This path is constant except on the subset $C'_\delta \subset C_\delta$ of points
satisfying $\max(|x_3| ,|x_1|) \leq 1+2\delta$. The situation is shown schematically in Figure 2. The one-form
\[ \alpha= x_2 x_1 dx_2 - x_4 x_1 dx_4  - x_4 x_2 dx_3. \]
satisfies $d\alpha = \zeta$. Moreover, $d(x_2^2+x_4^2)\wedge\alpha \wedge d\alpha = 4(x_2^2+x_4^2)dx_{1234}$,
which shows that $\alpha$ is a contact form on each positive level set of $(x_2^2+x_4^2)$.
Thus $D_\delta^* \alpha$ is contact on $C_\delta$.

Recall that $\Omega_t^\delta$ is the sum of terms $D_\delta^* \zeta$, $\delta\epsilon d( x_4 f_t^\delta dx_2)$,
and $\delta ^2 \epsilon f_t^\delta dx_{13}$.
We need primitives $\gamma_t^\delta$ for $f_t^\delta dx_{13}$. We choose $\gamma_t^\delta$ so that its
value at $x$ depends only on the values of $f_t^\delta$ along the path  $s \mapsto (x_1,s x_2,x_3,s x_4)$.
In particular, $\gamma^f_t $ is independent of $t$ at points of $C$ where $\max(|x_1|, |x_3|)  > 1+2\delta$.

Now consider the one-forms $D_\delta^* \alpha + \frac{\delta\epsilon}{2} (x_4 f_t^\delta dx_2+ \delta\gamma^\delta_t)$.
By decreasing $\epsilon$ (independently of $\delta$), we can ensure that all these forms are contact on $C_\delta$.

Finally, apply the stability lemma to the path of contact forms
$D_\delta^* \alpha + \frac{\delta\epsilon}{2} (x_4 f_t^\delta dx_2 + \delta \gamma^\delta_t)$ on $C_\delta$.
This path is constant outside the subset $C'_\delta$. The output, in the version described in Remark \ref{rephrasing},
is a new path of two-forms on $U_\delta$,
constant near $\partial U$, which can be plugged into our near-symplectic manifold near the chosen letter H.

This completes the proof of the theorem.

\section{The number of even circles}

\subsection{Oriented two-plane fields on three-manifolds}
The set $J(Y)$ of homotopy classes of smooth, oriented two-plane fields $\xi \subset TY$ over an oriented
three-manifold $Y$ is non-empty because $e(TY)=0$. A good way to describe its structure
(one which does not involve choosing a trivialisation of $TY$) is to exploit the classification of $\spinc$-structures.
We review this, following Kronheimer and Mrowka \cite{KM}.

On an oriented Riemannian three-manifold, a $\spinc$-structure $\tilde{\mathfrak{t}}$ is specified by a Hermitian
$\C^2$-bundle $\mathbb{S}\to Y$ and a map $\rho\colon T^*Y \to \mathfrak{su}(\mathbb{S})$
satisfying $\rho(\alpha)^2 = - |\alpha|^2 \id$. It should also be consistent with the orientation, which means that
$\rho(e_1)\rho(e_2) \rho(e_3)= \id_{\mathbb{S}}$ when $(e_1,e_2,e_3)$ is an oriented orthonormal basis for $T^*_yY$.

The isomorphism class of $\tilde{\mathfrak{t}}$ is an element $\mathfrak{t}$ of a freely transitive $H^2(Y;\Z)$-set $\spinc(Y)$. This
set is essentially independent of the metric.

\begin{Prop}
Oriented two-plane fields are in canonical bijection with isomorphism classes of pairs $(\tilde{\mathfrak{t}},\Psi)$,
where $\Psi\in \Gamma(Y,\mathbb{S})$ is a unit-length spinor.
\end{Prop}
\begin{pf}
The two-plane field associated with $(\tilde{\mathfrak{t}},\Psi)$ is $\ker (\alpha)$,
the kernel of the unique one-form $\alpha$ with $|\alpha|^2=1$ such that $\rho(\alpha)$ has $\C \Psi$ and $\Psi^\perp$
as $\ii$- and $-\ii$-eigenspaces. Conversely, an oriented two-plane field $\xi$ supports a hermitian structure,
unique up to homotopy, which induces a hermitian structure on the real vector bundle
$\mathbb{S}:= T^*Y\oplus \epsilon^1 = \xi^*\oplus (\xi^*)^\perp  \oplus \epsilon^1$.
Here $\xi^*\subset T^*Y$ is the image of $\xi$ under the metric isomorphism $TY \to T^*Y$, and $\epsilon^1$ is the
trivial real line bundle. The last summand has an obvious unit-length section $\Psi$.
The Clifford map $\rho$ is given, at each point $y\in Y$, by a map $L \oplus \R \to \mathfrak{su}(L\oplus \C)$, where $L$ is a
1-dimensional hermitian vector space; as such, it can be described using the quaternions $\mathbb{H}$.
Right multiplication by $j$ makes $\mathbb{H}$ a complex vector space, with complex subspaces
$\langle i,k\rangle$ and $\langle 1,j \rangle$.
Left multiplication gives a map $\langle i,k \rangle \oplus \langle j \rangle \to \mathfrak{su}(\mathbb{H})$, which gives the
model for $\rho$.
\end{pf}
It is straightforward to deduce the following results from this correspondence.
We refer to \cite{KM} for further explanation.

Frst, there is a coarse decomposition of $J(Y)$ by isomorphism classes of $\spinc$-structures,
slightly finer than that by Euler class:
\[  J(Y) = \bigcup_{\mathfrak{t}\in \spinc (Y)}{J(Y,\mathfrak{t})}.  \]
For a pair of two-plane fields $(\xi_0,\xi_1)$, representing $(\mathfrak{t}_0,\mathfrak{t}_1)$, the
difference $c_1(\mathfrak{t}_0) - c_1(\mathfrak{t}_1) \in H^2(Y;\Z)$ is the classical obstruction to homotopy of
$\xi_0$ and $\xi_1$ over the two-skeleton, which lies in $H^2(Y;\pi_2(S^2))=H^2(Y;\Z)$.

There is then a fine decomposition of $J(Y,\mathfrak{t})$: We assume $Y$
connected. Then each $J(Y,\mathfrak{t})$ is a transitive $\Z$-set: the action $j\mapsto j[n]$ of $n\in \Z$
is induced by an automorphism of $TY$, supported over a small ball in $Y$:
take a map $(B^3,\partial B^3)\to (\SO(3),1)$ of degree $2n$, and transfer it to $TY$ via an embedding
$B^3\hookrightarrow Y$.

The stabiliser of this action, a subset of $2\Z$, is the image of the homomorphism
\[ H^1(Y)\to \Z,\quad a\mapsto (a\cup  c_1(\mathfrak{t} )) [Y] .\]

\subsection{Almost complex structures} On an oriented Riemannian four-manifold these also have a spinorial interpretation
\cite{KM}:

\begin{Prop}
There is a canonical bijection between orthogonal, positively oriented almost complex structures and isomorphism classes of pairs
$(\tilde{\mathfrak{s}},\Phi)$ of $\spinc$-structure and unit-length positive spinor $\Phi\in \Gamma(\mathbb{S}^+)$.
\end{Prop}
The $\spinc$ structure $\tilde{\mathfrak{s}}$ is specified by hermitian $\C^2$-bundles $\mathbb{S}^+$, $\mathbb{S}^-$ with a
map $\rho\colon T^*X \to \Hom(\mathbb{S}^+,\mathbb{S}^-)$ such that $\rho^*(\alpha)\rho(\alpha)=|\alpha|^2\id_{\mathbb{S}^+}$.

We now depart from our review of Kronheimer and Mrowka's ideas. Take a compact, oriented four-manifold $X_0$, with non-empty boundary, and introduce
the set $J(X_0,\mathfrak{s})$ of homotopy classes of almost complex structures on $X_0$ underlying $\mathfrak{s}\in \spinc(X_0)$.
This set admits an action by $H_1(X_0,\partial X_0;\Z)$: for $a\in J(X_0,\mathfrak{s})$, and an oriented path
$\gamma$ with boundary in (and transverse to) $\partial X_0$, define $a[\gamma] $ by choosing a map
$ (B^3, \partial B^3)\to ( \mathrm{SO}(3),1)$ of degree 2, and using this to induce an automorphism
$\alpha$ of $TX_0$ supported in a tubular neighbourhood $B^3\times [0,1]$ of $\im(\gamma)$.
Choose a representative $I$ for $a$, and let $a [\gamma] $ be the class of $\alpha^*I $.
Note that this does represent $\mathfrak{s}$.

An almost complex structure preserves a unique field $[\xi]\in J(\partial
X_0)$ over its boundary. Thus there is a `restriction' map $J(X_0,\mathfrak{s})\to J(\partial X_0, \mathfrak{s}|_{\partial X_0})$.
The group actions are compatible:
$(a[\gamma])|_{\partial X_0}   = (a|_{\partial X_0}) [\partial \gamma ] $, where $[\partial \gamma]\in H_0(\partial X_0;\Z)$.

\begin{Lem} The $H_1(X_0,\partial X_0;\Z)$-action on $J(Y,\mathfrak{s})$ is transitive.
Its stabiliser is the image of the homomorphism $H_3(X_0,\Z)\to H_1(X_0;\Z)$, $x\mapsto c_1(\mathfrak{s})\cap x $.
\end{Lem}

\begin{pf}
We can represent orthogonal almost complex structures $J_0$, $J_1$ in the class $\mathfrak{s}$ by spinors
$\Psi_0$ and $\Psi_1=\theta \Psi _0 $, for a function  $\theta\colon X_0 \to S^1$. Homotopy of $\Psi_0$ and $\Psi_1$
is detected by the relative Euler class of $\mathbb{S}^+\times[0,1] \to X_0\times [0,1]$:
\begin{equation}
e(\mathbb{S}^+\times[0,1] ; \Psi_0,\Psi_1  )\in H_1(X_0\times[0,1], \partial X_0 \times [0,1] ;\Z).
\end{equation}
Indeed, given a cell decomposition of $X_0$ relative to $\partial X_0$,
$ e(\mathbb{S}^+\times[0,1] ; \Psi_0,\Psi_1  )$ is the obstruction to homotopy over the 3-skeleton;
but there is then no obstruction to extending the homotopy over the 4-skeleton.

Put $\delta(\Psi_0, \Psi_1) = i (e)$, where
\[i\colon  H_1(X_0\times[0,1], \partial X_0 \times [0,1] ;\Z)\to H_1(X_0,\partial X_0;\Z)\]
is the obvious isomorphism.
Thus $\delta(\Psi_0,\Psi_1)$  is represented by the projection to $X_0$ of the 1-cycle of zeros of a section of
$\mathbb{S}^+\times [0,1]$ extending the $\Psi_i$.

Now, if $a\in J(X,\mathfrak{s})$ is represented by $(\tilde{\mathfrak{s}},\Psi_0)$ then $a[\gamma]$ is represented
by a pair $(\tilde{\mathfrak{s}},\Psi_1)$ with $\delta(\Psi_0, \Psi_1)= [\gamma]$.
For, if $\alpha\in \aut(TX)$ induces the action by $\gamma$, then (using a standard model for $\alpha$)
the zeros of $(1-t)\Psi_0+ t\alpha^*\Psi_0$ occur only at $t=1/2$ and along $\gamma$.
Hence $J(X_0,\mathfrak{s})$ is a transitive $H_1(X_0,\partial X_0;\Z)$-set.

The function $\theta$ represents a class $[\theta]\in H^1(X_0,\partial X_0 ;\Z)$, and
\[  \delta(\Psi_0, \theta \Psi_0) = \pm c_1(\mathfrak{s})\cap  \mathrm{PD} [\theta ] \in H_1(X_0;\Z)\quad \text{(universal sign).} \]
This can be seen using the linear homotopy $(1-t) \Psi_0 + t \theta \Psi_0'$, where $\Psi'$ is a unit-length spinor
homotopic to $\Psi$ such that $\Psi_0\wedge \Psi_0'$ vanishes transversely.
Its zero-set is concentrated at $t=1/2$; there, it is the intersection of the zero-set of $\Psi_0\wedge \Psi_0'$
(representing $c_1(\mathbb{S}^+)$) and a level set of $\theta$. Thus the stabiliser is as claimed.
\end{pf}

\subsection{Application to near-symplectic geometry}

We can now read off from the discussion above that there is an isomorphism of $\Z$-sets $J(S^3)\cong \Z$, and that for $S^1\times S^2$ we have

\begin{Lem}
$J(S^1\times S^2) = \bigcup_{n\in \Z} J(S^1\times S^2, \mathfrak{t}_n )$, where $\mathfrak{t}_n\in \spinc(S^1\times S^2)$
is characterised by $\langle  c_1(\mathfrak{t}_n),[\{t\}\times  S^2]\rangle= 2n$. As $\Z$-sets,
\[ J(S^1\times S^2, \mathfrak{t}_n ) \cong \Z/(2n).\]
\end{Lem}

Choose disjoint embeddings $\bar{B}^4\hookrightarrow S^4$ and $S^1\times \bar{B}^3\hookrightarrow S^4$, and let $W$ be the closure of the complement of their images.

\begin{Lem} \label{plane fields}
We have
\[J(W)=\bigcup_{n\in \Z}{J(W,\mathfrak{s}_n)}, \]
where $\mathfrak{s}_n\in \spinc(W)$ is the unique element which restricts to $\mathfrak{t}_n$ on $S^1\times S^2$; and $H_1(W,\partial W;\Z) =\Z$, where $1\in \Z$ corresponds to a path $\gamma$ running from $S^3$ to $S^1\times S^2$.  Each $J(W,\mathfrak{s}_n)$ is acted upon transitively by $\Z$, and because $H^1(W;\Z)=0$ this action is free.

\end{Lem}

There are restriction maps $r_1 \colon J(W,\mathfrak{s}_n)\to J(S^3)$ and
$r_2  \colon J(W,\mathfrak{s}_n)\to  J(S^1\times S^2,\mathfrak{t}_n)$, satisfying
\[r_1(a[1]) =  (r_1(a))[-1] ,\quad r_2(a[1]) = (r_2(a))[1] .\]
To complete the description of the restriction maps it is necessary only
to give one pair $(r_1(a),r_2(a))$ for each $n$. The case which will be relevant to our purposes is $n=-1$.
The sets $J(S^3) $ and $J(S^1\times S^2,\mathfrak{t}_{-1})$ have distinguished basepoints,
represented respectively by the boundary of the standard complex structure on $B^4$,
and by an $S^1$-invariant 2-plane field. Referring to these basepoints, we have
\begin{Lem}\label{restr}
There is an $a$ with  $(r_1(a),r_2(a))=(0,0)\in \Z\times \Z/2$.
\end{Lem}

Indeed, this follows from the following calculation:

\begin{Lem}\label{Hopf calculation}
$\Theta_{\mathrm{ev}}$, restricted to $\{|x|=1\}$ and considered as a map $S^1\times S^2\to \R^3\setminus \{0\}$,
extends to a map  $ \bar{D}^2\times S^2 \to \R^3\setminus \{0\}$.
\end{Lem}
\begin{pf}
A tubular neighbourhood of a circle $S^1\subset S^4$ gives a decomposition
\[S^4 = (S^1\times \bar{D}^3) \cup_{S^1\times S^2} (\bar{D}^2\times S^2).\]
What follows is an explicit model for this decomposition.

Consider the function $f: \C^2\to [0,\infty)$, $f(z_1,z_2)=(|z_1|-2)^2+|z_2|^2$. The pair $(f^{-1}[0,1], f^{-1}(0))$ is
diffeomorphic to $(S^1\times \bar{D^3}, S^1\times\{0\})$, via the map
\[ \alpha : (z_1,z_2)\mapsto \left(|z_1|^{-1} z_1 ; |z_1|-2, \real(z_2),\imag(z_2)\right).\]
The pullback $\alpha^*\Theta_{\mathrm{ev}}$ is given by $(z_1,z_2)\mapsto (-2(|z_1|-2),z_2)$. Define
\[ F\colon  \C^2\setminus f^{-1}(0)\to \R^3\setminus \{0\};\quad (z_1,z_2)\mapsto f(z_1,z_2)^{-1/2}(-2(|z_1|-2),z_2).  \]
Then $F=\alpha^*\Theta_{\mathrm{ev}}$ on $f^{-1}(1)$. We consider the map $F/|F|$, restricted to a sphere
$S^3_R=\{|z_1|^2+|z_2|^2=R^2 \}$, $R\gg 0$. This map is not surjective ($(1,0,0)$ is not in its image).
But a map $S^3\to S^2$ which is not surjective has trivial Hopf invariant, and hence extends over the four-ball.
\end{pf}

\begin{pf}[Proof of \ref{Gompf}]

Let $n_0$ (resp. $n_1$) denote the number of even (resp. odd) circles of $\omega$.
Choose a Riemannian metric making $\omega$ self-dual. We proceed in five steps.

{\bf 1.}
We recall from \cite{HH} a corollary of Hirzebruch's signature theorem:
When one has disjointly embedded closed four-balls $B_1,\dots, B_N \subset X$, and a unit-length section
$\lambda$ of $\Lambda^+$ over $X\setminus \bigcup{\mathrm{int}(B^4)}$ with `Hopf invariants'
$h_1(\lambda),\dots h_N(\lambda)$, the relation
\[ (c^2 - 2e(X)-3\sigma(X))/4 = -\sum_{i=1}^N{h_i(\lambda)}  \]
holds. Here $c\in H^2(X;\Z)$ is the unique class which restricts to $c_1(X\setminus \bigcup B_i, \lambda)$,
over the complement of $\bigcup{B_i}$. The Hopf invariant $h_i(\lambda) \in [S^3,S^2] = \Z$ is the obstruction to
extending $\lambda$ over $B_i$. It could also be viewed as a relative Euler number $e(\mathbb{S}^+, \lambda)$
for the positive spinor bundle over the four-ball, with the non-vanishing spinor determined by $\lambda$ over the boundary.

Since $c$ is characteristic, $c^2 \equiv \sigma(X) \mod 8$, and so the relation reduces modulo 2 to
\[ \beta(X)\equiv \sum {h_i(\lambda)} \mod 2, \]
where we write $\beta(X)$ for the integer $1-b_1(X)+ b_2^+(X)$.

{\bf 2.} Let $X_\gamma$ denote the manifold obtained by surgery along an embedded, framed, oriented 1-submanifold
$\gamma\subset X$. This means that we excise an open tubular neighbourhood of $\gamma$ to obtain a manifold $X_0$; the framing gives (up to homotopy) an identification of $\partial X_0$ with a number of copies of $S^1\times S^2$. To each
boundary component we attach the complement of a standard copy of $B^3\times S^1\subset S^4$ (this is diffeomorphic
to $S^2\times B^2$). Then
\[ \beta(X_\gamma) = \beta(X)+1.\]
Indeed, $\beta=(e+\sigma)/2$. As a cobordism invariant, the signature $\sigma$ is unaffected by surgery, while a pair of
Mayer-Vietoris sequences yields
\[ e(X')- e(X) = e(S^2\times B^2)-e(B^3\times S^1)=2-0.\]

{\bf 3.}  We now carry out surgery along $Z_\omega$. Along each component there are two homotopically distinct framings,
since $\pi_1(\SO(3))=\Z/2$, but either will do here. Note that $\beta(X')=\beta(X)+n_0+n_1$.

{\bf 4.}  The result of this surgery is a manifold
\[ X'  = X_0 \cup \bigcup_{i=1}^N {X_i'}, \]
with standard pieces $X_i'$, and a non-vanishing self-dual two-form $\omega|X_0$ on $X_0$.
Let $I$ be an almost complex structure compatible with $\omega|X_0$.
It restricts to $\partial X_i'$ to give a class in $J(S^1\times S^2)$. This
class lies in $J(S^1 \times S^2,\mathfrak{t}_{-1})$, since as shown by Taubes \cite{Ta2},
$\langle c_1(\omega),[\{t\}\times S^2] \rangle = -2$ with our orientation conventions.
Because of our choice of framing, this class is the odd or even element of $J(S^1\times S^2, \mathfrak{t}_{-1})=\Z/2$
according to the parities of the original circles.

By Lemma \ref{restr}, the element $a \in J(S^1\times S^2, \mathfrak{t}_{-1})$ extends  over the standard manifold $W$ to
give an element $b\in J(S^3)=\Z$ if and only if $a= b\mod 2$.

{\bf 5.}
The previous steps together give (\ref{Gompf}): the form $\omega|X_0$ extends to a unit-length section of
$\Lambda^+ (X'\setminus\bigcup _{i=1}^{n_0+n_1}{B_i})$ (where $B_i$ is a small ball in $X_i'$).
The sum of the Hopf invariants is $n_1 \mod 2$, and we have $\beta(X)+n_0+n_1 \equiv n_1 \mod 2$.
\end{pf}


\begin{thebibliography}{99}
\bibitem{ADK}
D. Auroux, S. Donaldson and L. Katzarkov, \emph{Singular Lefschetz pencils}, Geom. Topol. 8 (2005), 1043--1114.

\bibitem{Cal}
E. Calabi, \emph{An intrinsic characterization of harmonic one-forms},  1969  Global Analysis (Papers in Honor of K. Kodaira),
101--117, Univ. Tokyo Press, Tokyo.

\bibitem{GK}
D. Gay and R. Kirby, \emph{Constructing symplectic forms on 4-manifolds which vanish on circles}, Geom. Topol. 8 (2004), 743--777.

\bibitem{HH}
F. Hirzebruch and H. Hopf, \emph{Felder von Fl\"achenelementen in 4-dimensionalen Mannigfaltigkeiten}, Math. Annalen 136 (1958), 156--172.

\bibitem{Hon}
K. Honda, \emph{Local properties of self-dual harmonic two-forms on a 4-manifold}, J. Reine Angew. Math. 577 (2004), 105--116.

\bibitem{Ho2}
K. Honda, \emph{Transversality theorems for harmonic forms}, Rocky Mountain J. Math. 34 (2004), no. 2, 629--664.

\bibitem{KM}
P. Kronheimer and T. Mrowka, \emph{Monopoles and contact structures}, Invent. Math. 130, 209--255 (1997).

\bibitem{Leb}
C. LeBrun, \emph{Yamabe constants and the perturbed Seiberg-Witten equations}, Comm. Anal. Geom. 5 (1997), no. 3, 535--553.

\bibitem{LS}
K. Luttinger and C. Simpson, \emph{A normal form for the birth/flight of closed selfdual 2-form degeneracies}, ETH preprint, unpublished, 1996.

\bibitem{Ta1}
C. Taubes, \emph{The geometry of the Seiberg-Witten invariants}, Proceedings of the International Congress of Mathematicians, Vol. II (Berlin, 1998). Doc. Math. 1998, Extra Vol. II, 493--504.

\bibitem{Ta2}
C. Taubes, \emph{Seiberg-Witten invariants and pseudo-holomorphic subvarieties
for self-dual, harmonic $2$-forms}, Geom. Topol. 3 (1999), 167--210.

\bibitem{Ta3}
C. Taubes, \emph{A proof of a theorem of Luttinger and Simpson about the number of vanishing circles of a near-symplectic form on a
4-dimensional manifold}, Math. Res. Lett.  13  (2006),  no. 4, 557--570.

\end{thebibliography}
\end{document}